\documentclass[11pt]{article}

\usepackage{amsmath,amsfonts,amsthm,amssymb,graphicx,tikz,tikz-cd,url,hyperref,cleveref}
\usepackage[all,2cell,ps]{xy}

\theoremstyle{plain}
\newtheorem{thm}{Theorem}[section]

\newtheorem{lem}[thm]{Lemma}
\newtheorem{prop}[thm]{Proposition}

\newtheorem{cor}[thm]{Corollary}

\newtheorem{qtn}[thm]{Question}

\theoremstyle{definition}

\newtheorem{ex}[thm]{Example}

\newtheorem{rem}[thm]{Remark}

\theoremstyle{remark}

\newcommand{\bbA}{\mathbb{A}}
\newcommand{\bbB}{\mathbb{B}}
\newcommand{\bbC}{\mathbb{C}}

\newcommand{\bbH}{\mathbb{H}}

\newcommand{\bbN}{\mathbb{N}}

\newcommand{\bbQ}{\mathbb{Q}}
\newcommand{\bbR}{\mathbb{R}}

\newcommand{\bbZ}{\mathbb{Z}}

\newcommand{\calC}{\mathcal{C}}

\newcommand{\calG}{\mathcal{G}}

\newcommand{\calI}{\mathcal{I}}

\newcommand{\calN}{\mathcal{N}}
\newcommand{\calO}{\mathcal{O}}
\newcommand{\calP}{\mathcal{P}}
\newcommand{\calQ}{\mathcal{Q}}

\newcommand{\frakg}{\mathfrak{g}}

\newcommand{\frakl}{\mathfrak{l}}

\newcommand{\al}{\alpha}
\newcommand{\gam}{\gamma}
\newcommand{\Gam}{\Gamma}

\newcommand{\Del}{\Delta}

\newcommand{\Lam}{\Lambda}
\newcommand{\sig}{\sigma}

\newcommand{\om}{\omega}

\DeclareMathOperator{\U}{U}

\DeclareMathOperator{\M}{M}
\DeclareMathOperator{\SL}{SL}
\DeclareMathOperator{\PSL}{PSL}
\DeclareMathOperator{\GL}{GL}

\DeclareMathOperator{\SO}{SO}

\DeclareMathOperator{\PU}{PU}
\DeclareMathOperator{\SU}{SU}

\DeclareMathOperator{\Sp}{Sp}

\DeclareMathOperator{\Id}{Id}

\newcommand{\bs}{\backslash}

\newcommand{\lra}{\longrightarrow}

\newcommand{\conj}{\overline}
\newcommand{\wh}{\widehat}
\newcommand{\wt}{\widetilde}

\newenvironment{pf}{\begin{proof}}{\end{proof}}

\newenvironment{enum}{\begin{enumerate}}{\end{enumerate}}

\allowdisplaybreaks

\makeatletter
\let\@@pmod\pmod
\DeclareRobustCommand{\pmod}{\@ifstar\@pmods\@@pmod}
\def\@pmods#1{\mkern4mu({\operator@font mod}\mkern 6mu#1)}
\makeatother

\usetikzlibrary{arrows}

\title{Residual finiteness and discrete subgroups of Lie groups \\ \bigskip\bigskip \small{\emph{Dedicated to the memory of Nicolas Bergeron}}}
\author{Matthew Stover \\ \small{Temple University}\\ \small{\textsf{mstover@temple.edu}}}
\date{\today}

\begin{document}

\maketitle

\begin{abstract}
Let $G$ be a real Lie group and $\Gam < G$ be a discrete subgroup of $G$. Is $\Gam$ residually finite? This paper describes known positive and negative results then poses some questions whose answers will lead to a fairly complete answer for lattices.
\end{abstract}

\section{Definitions and first examples}\label{sec:Intro}

A group $\Gam$ is \emph{residually finite} if for every nontrivial $\gam \in \Gam$ there is a homomorphism $\rho : \Gam \to F$ onto a finite group $F$ so that $\rho(\gam)$ is nontrivial. The applications of residual finiteness in algebra, number theory, geometry, and topology are abundant and significant. Just to give an example, there is \emph{Scott's criterion} \cite[Lem.\ 1.3]{Scott}:

\begin{itemize}
\item[] If $X$ is a Hausdorff space with fundamental group $\Gam$ and universal cover $\widetilde{X}$, then $\Gam$ is residually finite if and only if for every compact subset $C \subseteq \widetilde{X}$ there is a finite cover $Y \to X$ so that the natural map $\widetilde{X} \to Y$ restricts to an embedding on $C$.
\end{itemize}

\noindent
Thus residual finiteness can be used to promote certain immersions to embeddings in finite covers. The purpose of this paper is to present an account of what is known and what is expected regarding the residual finiteness of arguably the most important groups in all four of the aforementioned areas, namely discrete subgroups of Lie groups.

\begin{qtn}\label{qtn:Main1}
Let $G$ be a connected real Lie group. Which discrete subgroups $\Gam < G$ are residually finite?
\end{qtn}

It ends up that \Cref{qtn:Main1} is too broad. The remainder of this section provides some results and contrasting examples that motivate this survey's restriction to \emph{lattices}, i.e., discrete subgroups such that Haar measure on $G$ descends to a finite measure on the quotient $\Gam \bs G$. First, if $G$ admits a faithful linear representation, then Malcev famously showed that all finitely generated subgroups of $G$ are residually finite.

\begin{thm}[Malcev \cite{MalcevRF}]\label{thm:MalcevRF}
Let $\Gam$ be a finitely generated subgroup of $\GL_N(\bbC)$ for some $N \ge 1$. Then $\Gam$ is residually finite.
\end{thm}

The hypothesis that $\Gam$ is finitely generated is essential, since for example the trivial group is the only finite quotient of $\bbQ$. One even finds counterexamples among discrete subgroups. For example, while one can show that every discrete subgroup of $\PSL_2(\bbR)$ is residually finite, this does not even remain true for $\SL_2(\bbR)$.

\begin{ex}[Folklore]\label{ex:SL2RnotRF}
The linear Lie group $\PSL_2(\bbR)$ contains a discrete subgroup isomorphic to a countably infinite free product
\[
\Gam \cong \bbZ / 2 \ast \bbZ / 4 \ast \cdots \ast \bbZ / 2^j \ast \cdots
\]
of cyclic groups of order $2^j$. One constructs $\Gam$ using rotations $\sig_j \in \PSL_2(\bbR)$ of order $2^j$ fixing appropriate points $z_j$ in the upper half-plane $\bbH^2$, $j \in \bbN$. The preimage of each $\sig_j$ in $\SL_2(\bbR)$ is an element $\tau_j$ of order $2^{j+1}$ for which $\tau_j^{2^j} = -\Id$. Thus the preimage $\wh{\Gam}$ of $\Gam$ in $\SL_2(\bbR)$ fits into a central exact sequence
\[
1 \lra \langle \pm \Id \rangle \lra \wh{\Gam} \lra \Gam \lra 1
\]
induced by the cover $\SL_2(\bbR) \to \PSL_2(\bbR)$. If $\rho : \wh{\Gam} \to F$ is a finite quotient of $\wh{\Gam}$, then there is some $n_0$ so that $\rho(\tau_j)$ has order $2^{n_j}$ for some $0 \le n_j \le n_0$. As soon as $n_0 < j$, one obtains
\[
\rho(-\Id) = \rho\!\left(\tau_j^{2^j}\right) \!= \rho\!\left(\tau_j^{2^{n_j}}\right)^{2^{j-n_j}} \!= 1
\]
and therefore $-\Id$ maps trivially to every finite quotient of $\wh{\Gam}$.
\end{ex}

If $\Gam < \SL_2(\bbR)$ is discrete and torsion-free, then it is isomorphic to its image in $\PSL_2(\bbR)$ under the covering projection and thus is residually finite. One might hope based on this that torsion is essential for producing examples like \Cref{ex:SL2RnotRF} for more general Lie groups, but one only needs to replace $\bbR$ with $\bbC$ to eliminate that possibility.

\begin{ex}[Cremaschi and Souto \cite{CremaschiSouto}]
There is a discrete, torsion-free subgroup $\Gam < \PSL_2(\bbC)$ that is not residually finite.
\end{ex}

As indicated by \Cref{ex:SL2RnotRF}, residual finiteness can behave badly under finite coverings of Lie groups. This phenomenon is not restricted to infinitely generated monstrosities. The symplectic group $\Sp_{2 g}(\bbR) < \SL_{2g}(\bbR)$ is homotopy equivalent to the unitary group $\U(g)$, which has fundamental group $\bbZ$. Thus $\Sp_{2 g}(\bbR)$ has a connected $n$-fold cyclic cover for all $n \in \bbN$; these covers are all nonlinear real Lie groups. The following theorem of Deligne will be put in greater context in \Cref{sec:ss}.

\begin{thm}[Deligne \cite{DeligneSp}]\label{thm:Deligne}
Let $\Gam = \Sp_{2 g}(\bbZ)$ and $G_n \to \Sp_{2 g}(\bbR)$ be the unique finite, connected cover of degree $n$. If $g \ge 2$ and $n \ge 3$, then the preimage $\Gam_n$ of $\Gam$ in $G_n$ is a $\bbZ / n$ central extension of $\Gam$ that is not residually finite.
\end{thm}

Since $\Sp_{2 g}(\bbZ)$ is finitely generated, the linearity assumption in \Cref{thm:MalcevRF} is just as essential as finite generation of $\Gam$. In fact, Deligne gives a precise measure of the failure of residual finiteness. If $Z \cong \bbZ / n \trianglelefteq \Gam_n$ is the center, then the intersection of the kernels of all homomorphisms from $\Gam_n$ onto finite groups is $2 Z$. Thus \Cref{qtn:Main1} is even subtle for lattices. Indeed, $\Sp_{2 g}(\bbZ)$ is a lattice in $\Sp_{2 g}(\bbR)$ and its preimage $\Gam_n$ in the $n$-fold cyclic cover $G_n$ of $\Sp_{2 g}(\bbR)$ is a lattice in $G_n$. Having now seen that residual finiteness is even delicate for lattices in Lie groups, the following questions are the primary focus of this article.

\begin{qtn}\label{qtn:Main2}
For which connected real Lie groups $G$ are all lattices in $G$ residually finite?
\end{qtn}

\begin{qtn}\label{qtn:Main3}
Which connected real Lie groups contain a lattice that is not residually finite?
\end{qtn}

Both questions are quickly answered in the linear and solvable cases.

\begin{thm}\label{thm:LinearSolvable}
Let $G$ be a connected real Lie group that is either solvable or that admits a faithful linear representation. Then every lattice in $G$ is residually finite.
\end{thm}

\begin{pf}
If $G$ is linear, then every lattice in $G$ is finitely generated (see \cite[Cor.\ 1.6]{GelanderSlutsky} for a unified proof), hence lattices are residually finite by \Cref{thm:MalcevRF}. If $G$ is solvable, then Mostow showed that every lattice in $G$ is polycyclic \cite[Cor.\ 4.2]{Mostow} and Hirsch proved that polycyclic groups are residually finite \cite[Thm.\ 3]{Hirsch}. This proves the theorem.
\end{pf}

\begin{rem}\label{rem:Solvable}
Baumslag showed that not all finitely generated solvable groups are residually finite \cite{Baumslag}. These groups evidently are not isomorphic to a lattice in a solvable Lie group.
\end{rem}

\noindent
The remainder of the paper is organized as follows:
\begin{itemize}

\item[$\star$] Basic facts about residual finiteness and Lie groups are contained in \Cref{sec:Basics}.

\item[$\star$] The semisimple case, which is arguably the most varied and interesting, is considered in \Cref{sec:ss}. See Remarks \ref{rem:SArithmetic} and \ref{rem:SArithmetic2} for comments on extending the results described in this paper to the $S$-arithmetic setting of lattices in products of semisimple Lie groups over local fields.

\item[$\star$] The general case is studied in \Cref{sec:Gen}.

\end{itemize}

\subsubsection*{Acknowledgments} My paper \cite{StoverToledo2} with Domingo Toledo on residual finiteness of lattices in covers of $\PU(n,1)$ benefited immensely from correspondence with Nicolas Bergeron, who very sadly passed away in February 2024. This paper is dedicated to his memory. As everyone that had the pleasure to interact with him knows, his impact on mathematics goes far beyond his truly exceptional body of work, and he will be sorely missed by the community.

\medskip

I also thank Domingo Toledo for discussions, collaboration, and insight related to this paper and Richard Hill for related conversations. Thanks are also due to Daniel Groves for a conversation about relatively hyperbolic Dehn filling. Part of this paper was written while I was a Professeur Invit\'e at Universit\'e Grenoble-Alpes, and I thank them for their fantastic hospitality and working environment. Finally, this paper is loosely based on a lecture I gave at the workshop \emph{Zariski dense subgroups, number theory and geometric applications} at ICTP, Bangalore in January 2024. I want to thank the organizers for their persistence in realizing the workshop after the original incarnation was canceled by the pandemic and for the wonderful conference that resulted from their efforts. This material is based upon work supported by Grant Number DMS-2203555 from the National Science Foundation.

\section{Basic facts}\label{sec:Basics}

This section collects some basic facts used throughout this paper: \S\ref{ssec:RF} on residual finiteness and \S \ref{ssec:Lie} on Lie groups.

\subsection{Residual finiteness}\label{ssec:RF}

Recall that a group $\Gam$ is \emph{residually finite} if for all nontrivial $\gam \in \Gam$ there is a homomorphism $\rho$ onto a finite group $F$ so that $\rho(\gam)$ remains nontrivial. The next collection of facts are standard, and the uninitiated reader should prove them as warm-up exercises for what follows.

\begin{lem}\label{lem:FIRF}
Suppose that $\Gam$ is a group and $\Lam < \Gam$ is a finite index subgroup. Then $\Gam$ is residually finite if and only if $\Lam$ is residually finite.
\end{lem}

Two subgroups $\Gam, \Lam$ of a group $G$ are (widely) commensurable if there is a $g \in G$ so that $g \Lam g^{-1} \cap \Gam$ is finite index in $g \Lam g^{-1}$ and $\Gam$. Then one has:

\begin{cor}\label{cor:CommRF}
If $\Gam$ and $\Lam$ are commensurable subgroups of a group $G$, then $\Gam$ is residually finite if and only if $\Lam$ is residually finite.
\end{cor}

\begin{lem}[Malcev \cite{MalcevSemi}]\label{lem:SemiDirect}
Suppose that $\Lam = \Del \rtimes \Gam$ is a semidirect product of groups, with $\Del$ finitely generated. Then $\Lam$ is residually finite if and only if $\Del$ and $\Gam$ are residually finite.
\end{lem}

\begin{lem}\label{lem:FiniteInject}
Suppose that $\Gam$ is a residually finite group. For any finite subset $X \subseteq \Gam$ there is a homomorphism $\rho : \Gam \to F$ onto a finite group $F$ for which the restriction of $\rho$ to $X$ is injective.
\end{lem}

A important consequence of Lemma \ref{lem:FiniteInject} is the following corollary.

\begin{cor}[cf.\ Lem.\ 2.5 \cite{StoverToledo1}]\label{cor:JustCenter}
Suppose that $\Gam$ is a group for which there is a central extension
\[
1 \lra Z \lra \Gam \lra \Gam_0 \lra 1
\]
with $Z$ finite and $\Gam_0$ residually finite. Then $\Gam$ is residually finite if and only if there is a homomorphism $\rho : \Gam \to F$ onto a finite group $F$ so that $\rho$ is injective on $Z$.
\end{cor}

\begin{pf}
If $\Gam$ is residually finite, then the desired conclusion follows from Lemma \ref{lem:FiniteInject}. Suppose the converse. If $\gam \in \Gam$ has nontrivial image in $\Gam_0$, then residual finiteness of $\Gam_0$ provides a finite quotient of $\Gam_0$, which is thus a finite quotient of $\Gam$, where the image of $\gam$ is nontrivial. Thus it suffices to consider the case where $\gam \in Z$ is nontrivial, where the desired finite quotient exists by hypothesis.
\end{pf}

This section concludes with some remarks on using nilpotent quotients to prove residual finiteness. If a central extension
\[
1 \lra \bbZ \lra \wt{\Gam} \lra \Gam \lra 1
\]
is residually finite, then the reductions modulo $d \in \bbN$
\[
1 \lra \bbZ / d \lra \Gam_d \lra \Gam \lra 1
\]
are residually finite for \emph{some} infinite collection of $d$. The use of nilpotent quotients to prove residual finiteness of $\wt{\Gam}$ allows for even stronger conclusions. Residual finiteness of finitely generated nilpotent groups is particularly robust, thus nilpotent quotients of arbitrary finitely generated groups can be used to prove stronger residual finiteness statements. For example, the following result will play an important role later in this paper.

\begin{lem}[Lem.\ 2.6 \cite{StoverToledo1}]\label{lem:2StepQuo}
Suppose that $\Gam$ is a finitely generated residually finite group and that
\[
1 \lra \bbZ \lra \wt{\Gam} \lra \Gam \lra 1
\]
is a central exact sequence associated with $\phi \in H^2(\Gam, \bbZ)$. For $d \in \bbN$, let $\phi_d \in H^2(\Gam, \bbZ/d)$ be the reduction of $\phi$ modulo $d$ and
\[
1 \lra \bbZ / d \lra \Gam_d \lra \Gam \lra 1
\]
be the associated central exact sequence given by reducing the kernel $\langle \sig \rangle$ of $\wt{\Gam} \to \Gam$ modulo $d$. Suppose that there is a nilpotent quotient $\calN$ of $\wt{\Gam}$ that is injective on $\langle \sig \rangle$. Then $\Gam_d$ is residually finite for all $d \in \bbN$.
\end{lem}

The reader interested in the basic idea behind the proof should consider the case where $\calN$ is the Heisenberg group of $3 \times 3$ strictly upper-triangular integer matrices and where $\sig$ maps to the element
\[
\begin{pmatrix} 1 & 0 & c \\ 0 & 1 & 0 \\ 0 & 0 & 1 \end{pmatrix}
\]
for some nonzero $c \in \bbZ$. In a certain sense, the proof for that case (virtually) generalizes to arbitrary nilpotent groups by a famous theorem of Malcev \cite{MalcevNilpotent}, since all finitely generated nilpotent groups virtually embed in the subgroup of upper-triangular matrices in $\GL_N(\bbR)$ for some $N$.

\subsection{Lie groups}\label{ssec:Lie}

This section covers some basics on structure theory of Lie groups, along with some general consequences for residual finiteness of finitely generated subgroups.

\begin{thm}\label{thm:Liepi1}
Let $G$ be a connected real Lie group with universal cover $\wt{G}$. Then $\wt{G}$ is a Lie group and there is a natural homomorphism
\[
c : \pi_1(G) \to Z(\wt{G})
\]
to the center of $\wt{G}$ so that $\wt{G} \to \wt{G} / c(\pi_1(G)) \cong G$ is the universal covering.
\end{thm}

In particular, notice that the universal covering map is a homomorphism and that $\pi_1(G)$ is abelian, hence all connected covers of $G$ are regular and are naturally central extensions of $G$. From this point forward, the embedding of $\pi_1(G)$ into $Z(\wt{G})$ will be fixed and the homomorphism $c$ will be suppressed. \Cref{thm:Liepi1} has the following critical consequence for this paper.

\begin{cor}\label{cor:CentralDeck}
Let $G$ be a connected real Lie group and $G^\prime \to G$ be a connected cover with covering group $Z^\prime$. If $\Gam < G$ is a group and $\Gam^\prime$ is the preimage of $\Gam$ in $G^\prime$ under the covering, then $\Gam^\prime$ is a central extension
\[
1 \lra Z^\prime \lra \Gam^\prime \lra \Gam \lra 1
\]
of $\Gam$ by the abelian group $Z^\prime$.
\end{cor}

The next proposition reduces much of the work in this paper to understanding discrete subgroups of simply connected Lie groups with finite center.

\begin{prop}\label{prop:FiniteCoverRF}
Suppose that $G$ is a connected real Lie group and $\wh{G} \to G$ is a connected finite cover. If $\Gam < G$ has preimage in $\wh{G}$ that is residually finite, then the preimage of $\Gam$ in any connected cover of $G$ intermediate to $\wh{G}$ is residually finite.
\end{prop}

\begin{pf}
Let $Z < \wh{G}$ denote the deck group of the cover $\wh{G} \to G$, and consider the preimage $\wh{\Gam} < \wh{G}$ and a cover $G^\prime \to G$ intermediate to $\wh{G}$ with group $Z^\prime$. This induces a natural homomorphism of exact sequences
\[
\begin{tikzcd}
1 \arrow[r] & Z \arrow[r] \arrow[d, twoheadrightarrow] & \wh{\Gam} \arrow[r] \arrow[d, twoheadrightarrow] & \Gam \arrow[r] \arrow[d, "\mathrm{id}"] & 1 \\
1 \arrow[r] & Z^\prime \arrow[r] & \Gam^\prime \arrow[r] & \Gam \arrow[r] & 1
\end{tikzcd}
\]
with $\Gam^\prime$ the preimage of $\Gam$ in $G^\prime$.

Suppose that $\wh{\Gam}$ is residually finite. Then there is a homomorphism ${\rho : \wh{\Gam} \to F}$ onto a finite group $F$ so that $\rho$ is injective on $Z$ by \Cref{lem:FiniteInject}. If $\Lam$ is the kernel of $\rho$, then $\Lam \cap Z = \{1\}$. In particular, $\Lam$ intersects the kernel of the homomorphism $\wh{\Gam} \to \Gam^\prime$ in the identity, and thus it projects isomorphically to its image in $\Gam^\prime$. In other words, $\wh{\Gam}$ and $\Gam^\prime$ share an isomorphic subgroup of finite index. Applying \Cref{lem:FIRF} twice, one concludes that $\Lam$ and then $\Gam^\prime$ are both residually finite. This proves the proposition.
\end{pf}

This renders the case of linear universal cover relatively uninteresting for the purpose of this paper.

\begin{cor}\label{cor:LinearUC}
Suppose that $G$ is a connected real Lie group with finite fundamental group and linear universal cover. Then any finitely generated subgroup of $G$ is residually finite.
\end{cor}

Now, let $G$ be a semisimple Lie group. Then $G$ has a maximal compact subgroup $K$. The following famous theorem of \'{E}.\ Cartan reduces the problem of computing fundamental groups of semisimple Lie groups to that of their compact subgroups.

\begin{thm}[See Thm.\ VI.2.2 \cite{Helgason}]\label{thm:Cartan}
Let $G$ be a connected semisimple Lie group with maximal compact subgroup $K$. Then there is an $N \ge 1$ so that $G$ is diffeomorphic to $\bbR^N \times K$. In particular, the inclusion of $K$ into $G$ induces an isomorphism $\pi_1(G) \cong \pi_1(K)$.
\end{thm}

The remainder of this section discusses universal covers of (almost) simple Lie groups. Some preliminary definitions are required. Let $G$ be a real Lie group. Then $G$ is \emph{almost simple} if its center $Z(G)$ is finite and $G / Z(G)$ is simple in the usual sense of having no nontrivial proper normal subgroups. When $G$ admits an embedding into $\GL_N(\bbR)$ under which it is the vanishing set of a collection of polynomial equations on $\M_N(\bbR) \cong \bbR^{N^2}$, then $G$ is a \emph{linear algebraic group} over $\bbR$. More precisely, there is a linear algebraic group $\calG$ defined over $\bbR$ so that $G \cong \calG(\bbR)$. For a linear algebraic group, it makes sense to take its complexification $G(\bbC) = \calG(\bbC)$, which (as an abstract group) is independent of the embedding of $G$ in $\GL_N(\bbR)$.

When $G$ is a linear algebraic group over $\bbR$, it is \emph{absolutely almost simple} if $G$ is almost simple and its complex points $G(\bbC)$ remain almost simple. Similarly, $G$ is \emph{absolutely simply connected} if $G(\bbC)$ is simply connected\footnote{This is slightly nonstandard. In the literature one often says `absolutely almost simple, simply connected algebraic group' to mean simply connected in the algebraic sense. To be crystal clear about which simple connectivity is meant at a given point in this paper, the word absolutely is adjoined.} ; in the sense of algebraic geometry, the fundamental group of $G$ is the fundamental group over the algebraic closure. The defect between the topological and algebraic definitions of simple connectivity is determined by the following widely-known result.

\begin{lem}\label{lem:UCLinear}
Suppose that $G$ is a simple Lie group with universal cover $\wt{G}$. Then $\wt{G}$ is linear if and only if it is isomorphic to the (real points of the) unique absolutely almost simple, absolutely simply connected linear algebraic group $\calG$ over $\bbR$ that is locally isomorphic to $G$.
\end{lem}

\begin{pf}[Sketch of the proof]
Suppose that $\rho : \wt{G} \to \GL_N(\bbR)$ is a faithful linear representation. This determines a faithful representation of the Lie algebra $\frakg$ of $G$ into $\frakg\frakl_N(\bbR)$ whose complexification determines a representation $\rho_\bbC$ of the simply connected group $\calG(\bbC)$ to $\GL_N(\bbC)$ with central (hence finite) kernel. However then
\begin{align*}
\rho(\wt{G}) &= \rho_\bbC(\calG(\bbC)) \cap \GL_N(\bbR) \\
&= \rho_\bbC(\calG(\bbR))
\end{align*}
which implies that $\rho$ factors through the covering $\wt{G} \to \calG(\bbR)$. Since $\rho$ is faithful, it follows that $\wt{G} \cong \calG(\bbR)$.
\end{pf}

Finally, see the Appendix for tables enumerating the fundamental groups of simple Lie groups and other related facts like whether or not the universal cover is linear.

\section{Semisimple groups}\label{sec:ss}

This section describes what is known about and expected regarding residual finiteness of lattices in real semisimple Lie groups. As with much of the theory of lattices in semisimple groups, the discussion breaks naturally into two cases depending on the real rank of the ambient group. In the higher rank (i.e., real rank at least two) setting, assuming some widely-believed conjectures that are proved in many important cases, one obtains a very precise understanding of when irreducible lattices in arbitrary semisimple Lie groups are residually finite and when they are not. This is described in \S \ref{ssec:Higher}. The real rank one setting, where much is known but some concrete cases remain mysterious, is discussed in \S \ref{ssec:Rank1}.

\subsection{Higher rank}\label{ssec:Higher}

The following very basic example, first found by Millson \cite{Millson} almost simultaneous with Deligne's work \cite{DeligneSp}, hopefully hammers home that lattices that are not residually finite lurk very close to some of the most important (and obviously residually finite) lattices.

\begin{ex}[Millson \cite{Millson}]
Fix $n \ge 3$ and let $\calO$ be the ring of integers in the real quadratic field $\bbQ(\sqrt{p})$, where $p$ is a rational prime congruent to $7$ modulo $8$. If $\Gam = \SL_n(\calO)$ and $\wt{\Gam}$ is the preimage of $\Gam$ in the universal cover of $\SL_n(\bbR)$, then $\wt{\Gam}$ is not residually finite.
\end{ex}

Recall that $\SL_n(\bbR)$ is homotopy equivalent to $\SO(n)$, hence its fundamental group is $\bbZ / 2$ for $n \ge 3$. Also, note that $\Gam$ is actually a lattice in $\SL_n(\bbR) \times \SL_n(\bbR)$, not $\SL_n(\bbR)$. Nevertheless, any real embedding of $\bbQ(\sqrt{p})$ embeds $\Gam$ in $\SL_n(\bbR)$ and then $\wt{\Gam}$ is a central extension
\[
1 \lra \bbZ / 2 \lra \wt{\Gam} \lra \Gam \lra 1
\]
that is not residually finite, and the only obstruction is the lone nontrivial element of the central $\bbZ / 2$.

\medskip

Turning to the general case, if $\Gam$ is an irreducible lattice in a higher rank real Lie group $G$ with finite center $Z(G)$ and adjoint group $\conj{G} = G / Z(G)$, then $\Gam$ is arithmetic by famous work of Margulis \cite[p.\ 4]{Margulis}. Precisely, this means that there is an absolutely almost simple, absolutely simply connected linear algebraic group $\calG$ defined over a number field $k$ and a surjection $p$ from $\calG(k \otimes \bbR)$ onto $\conj{G}$ with compact kernel so that the image of $\Gam$ in $\conj{G}$ is commensurable with the image of $p(\calG(\calO_k))$, where $\calO_k$ is the ring of integers of $k$. Since $G$ is real, $k$ is a totally real number field.

Arithmetic properties of the algebraic group $\calG$ have a profound effect on residual properties of $\Gam$. The group $\calG(k)$ inherits two topologies from $\calG(\calO_k)$:
\begin{enum}

\item The topology induced by the \emph{profinite topology} on $\calG(\calO_k)$, where the basic open neighborhoods of $\Id \in \calG(k)$ are the finite index subgroups of $\calG(\calO_k)$.

\item The topology induced by the \emph{congruence topology} on $\calG(\calO_k)$, where now the basic open neighborhoods of $\Id \in \calG(k)$ are the finite index subgroups of $\calG(\calO_k)$ that contain the kernel of the reduction homomorphism $\calG(\calO_k) \to \calG(\calO_k / \calI)$ modulo a nonzero ideal $\calI \subseteq \calO_k$.

\end{enum}
Let $\wh{\calG}(k)$ denote the completion of $\calG(k)$ for the profinite topology and $\calG^\prime(k)$ denote its completion with respect to the congruence topology. Since congruence subgroups are finite index subgroups, there is an induced exact sequence
\[
1 \lra \calC(k) \lra \wh{\calG}(k) \lra \calG^\prime(k) \lra 1
\]
of continuous homomorphisms of profinite groups with kernel $\calC(k)$ that is called the \emph{congruence kernel}.

Then $\calG$ is said to have the \emph{congruence subgroup property} if $\calC(k)$ is central in $\wh{\calG}(k)$. It is a widely-famous conjecture, attributed to Serre \cite{Serre}\footnote{Here Serre does not conjecture, merely asks, whether the congruence subgroup property holds in higher-rank; see the footnote on p.\ 1.}, that $\calG$ always has the congruence subgroup property in higher rank. Moreover, this conjecture is known to hold in a great many cases; for instance, see \cite[\S 5]{PrasadRapinchukBMS} for a thorough account of where the conjecture is proved.

Let $\bbA_k$ be the adele ring of $k$. Restriction defines a homomorphism
\[
H^2(\calG(\bbA_k)) \lra H^2(\calG(k))
\]
whose kernel $\M(\emptyset, \calG)$ is known as the \emph{absolute metaplectic kernel}. The following, written closely following the introduction to \cite{PrasadRapinchuk}, was established by Deligne in the process of proving Theorem \ref{thm:Deligne}.

\begin{thm}[Deligne \cite{DeligneSp}]
Let $\calG$ be an absolutely almost simple, absolutely simply connected linear algebraic group over a number field $k$ that has the congruence subgroup property. Then no arithmetic subgroup in a cover of $\calG(k \otimes \bbR)$ of degree larger than the order of $\M(\emptyset, \calG)$ is residually finite.
\end{thm}

The absolute metaplectic kernel was completely determined\footnote{If $\calG$ is of type ${}^2\mathrm{A}_r$ arising from a noncommutative division algebra over a quadratic extension of $k$, one must assume Conjecture (U) from \cite[\S 2]{PrasadRapinchuk}. This conjecture presently remains open.} by Prasad and Rapinchuk \cite[Main Thm.]{PrasadRapinchuk}. For this paper, it suffices to say that $\M(\emptyset, \calG)$ is isomorphic to the dual $\wh{\mu}(k)$ of the group $\mu(k)$ of roots of unity of $k$. Since $k$ is totally real, $\mu(k) = \{\pm 1\}$. This implies that there are many lattices in higher rank Lie groups that are not residually finite. For example, the following result can be used in conjunction with Tables \ref{tb:Classicalpi1} and \ref{tb:Exceptionalpi1} to construct many examples.

\begin{thm}\label{thm:CSP-RF}
Suppose that $G$ is a connected semisimple Lie group with real rank at least two and $\Gam < G$ is an irreducible lattice with associated absolutely almost simple, absolutely simply connected algebraic group $\calG$ defined over the totally real number field $k$. Let $G_0$ denote the connected, absolutely almost simple, absolutely simply connected real linear algebraic group locally isomorphic to $G$, and suppose that $G$ is a covering of $G_0$ with degree greater than $\#\mathrm{M}(\emptyset, \calG) \le 2$. If $\calG$ has the congruence subgroup property, then $\Gam$ is not residually finite.
\end{thm}

Note that lattices in nonlinear groups $G$ covered by $G_0$ may be residually finite by Proposition \ref{prop:FiniteCoverRF}. Indeed, $G_0$ is linear and has finite center. To point out one notable special case in the opposite direction, lattices in groups with infinite fundamental group should never be residually finite (cf.\ \cite{DeligneSp, RaghunathanCentral}).

\begin{cor}\label{cor:CSP-infinite}
Let $G$ be a simple Lie group of higher rank with infinite fundamental group, and assume the congruence subgroup property for lattices in the absolutely almost simple, absolutely simply connected real linear algebraic group whose real points are locally isomorphic to $G$. Then no lattice in $G$ is residually finite.
\end{cor}

Thus the story in higher rank is fairly well-understood, modulo the congruence subgroup problem. In fact, the following remark indicates that this understanding is more general.

\begin{rem}\label{rem:SArithmetic}
While this paper is focused on real Lie groups, the considerations in this section should also give a complete understanding for irreducible $S$-arithmetic lattices in higher rank semisimple groups that are products of groups over local fields of characteristic zero. Note that even allowing the field $k$ to have complex places makes $\mu(k)$ and hence $\mathrm{M}(\emptyset, \calG)$ potentially more complicated. Interestingly, fixing a product of simple groups but changing the field $k$ could possibly influence the degree of the covering to which one must pass to find a lattice that is not residually finite. A motivated reader should be able to assemble precise statements and expectations from the references in this section.
\end{rem}

\subsection{Rank one}\label{ssec:Rank1}

The real Lie groups of real rank one are locally isomorphic to $\SO^\circ(n, 1)$, $\SU(n,1)$, $\Sp(n,1)$, or the adjoint group $\mathrm{F}_4^{(-20)}$. This section is broken up by group, which roughly divides into complete, partial, and very limited knowledge. A final subsection briefly notes the connection between this section and one of the most important unresolved problems in geometric group theory.

\subsubsection{Full understanding: $\SO^\circ(n,1)$ and $\Sp(n,1)$}

For $n \ge 3$ the universal cover of $\SO^\circ(n,1)$ is $\mathrm{Spin}(n,1)$, which is a linear finite cover, and $\Sp(n,1)$ is simply connected and linear for all $n \ge 2$. It follows from Proposition \ref{prop:FiniteCoverRF} that these examples are discarded with relative ease.

\begin{lem}\label{lem:EasyRank1}
Let $G$ be a connected Lie group locally isomorphic to $\SO^\circ(n,1)$ for $n \ge 3$ or $\Sp(n,1)$ for $n \ge 2$. If $\Gam < G$ is a lattice, then $\Gam$ is residually finite.
\end{lem}

The case $\PSL_2(\bbR) \cong \SO^\circ(2,1)$ is very special, first because it is arguably the most important simple Lie group, but also in the context of this paper since its fundamental group is $\bbZ$. This case is quite classical and there are arguments employing any of the many points of view that intersect in and around Riemann surfaces (see \cite[\S IV.48]{delaHarpe} for an account). One argument is sketched here, since it gives some preparation for what comes later in this section where an analogous argument will be given for certain other lattices.

\begin{thm}\label{thm:SL2RF}
Let $G$ be $\PSL_2(\bbR)$ and $G^\prime \to G$ be any connected covering. If $\Gam^\prime < G^\prime$ is a lattice, then $\Gam^\prime$ is residually finite.
\end{thm}

\begin{pf}[Sketch of proof]
Let $\Gam$ be the image of $\Gam^\prime$ in $\PSL_2(\bbR)$. By Corollary \ref{cor:CommRF} it suffices to consider the case where $\Gam$ is the fundamental group of a hyperbolic surface of finite area. Consider the case where $G^\prime = \wt{G}$ is the universal cover. Then the standard identification of $\wt{G}$ with the unit tangent bundle to $\bbH^2$ identifies the preimage $\wt{\Gam}$ of $\Gam$ in $\wt{G}$ with the fundamental group of the unit tangent bundle to $\Gam \bs \bbH^2$.

If $\Gam \bs \bbH^2$ is noncompact, then $\wt{\Gam} \cong \bbZ \times \Gam$, which is certainly residually finite, as is the preimage $\bbZ / d \times \Gam$ of $\Gam$ in any $d$-fold cover of $G$. It remains to consider the case where $\Gam \bs \bbH^2$ is closed with genus $g \ge 2$. The unit tangent bundle is the circle bundle with Euler number $2(1-g)$, which implies that $\wt{\Gam}$ has presentation
\[
\wt{\Gam} =\! \left\langle a_1, b_1, \dots, a_g, b_g, z\ \Big|\ \prod_{j = 1}^g [a_j, b_j] = z^{2(1-g)},\, z \textrm{ central } \right\rangle,
\]
which has an associated exact sequence
\[
1 \lra \langle z \rangle \lra \wt{\Gam} \lra \Gam \lra 1
\]
compatible with $\wt{G} \to G$. Taking the quotient of $\wt{\Gam}$ by the normal subgroup generated by all generators but $a_1, b_1, z$ gives a surjective homomorphism onto the group
\[
\calN =\! \left\langle \al, \beta, \xi\ \Big|\ [\al, \beta] = \xi^{2(1-g)} \right\rangle,
\]
which is a two-step nilpotent group. In fact, $\calN$ is isomorphic to a finite index subgroup of the integer Heisenberg group of strictly upper-triangular $3 \times 3$ integer matrices.

Note that the surjection $\wt{\Gam} \to \calN$ is injective on the center $\langle z \rangle$ of $\wt{\Gam}$. Since $\Gam < \PSL_2(\bbR)$ is residually finite, Lemma \ref{lem:2StepQuo} implies that $\wt{\Gam}$ and its reductions modulo $d \in \bbN$ are residually finite for all $d$. The reduction modulo $d$ of the center is the preimage of $\Gam$ in the $d$-fold cyclic cover of $\PSL_2(\bbR)$, so this proves the theorem.
\end{pf}

\begin{rem}\label{rem:SArithmetic2}
To complement Remark \ref{rem:SArithmetic2} and complete the picture for lattices in products of semisimple groups over local fields of characteristic zero, lattices in adjoint rank one groups over nonarchimedean fields of characteristic zero are virtually free, since the associated Bruhat--Tits building is a tree \cite[\S 2.7]{Tits}. Any extension of a free group by a finitely generated residually finite group is itself residually finite \cite[Cor.\ 7.11]{Reid}. Therefore lattices in rank one groups over $p$-adic fields are always residually finite.
\end{rem}

\subsubsection{Positive results for $\PU(n,1)$}

Note that $\PSL_2(\bbR)$ is also locally isomorphic to $\PU(1,1)$. For the purposes of this paper it ends up that $\PU(n,1)$ is more closely related to $\PSL_2(\bbR)$ than $\SO^\circ(n,1)$, in particular since it has fundamental group $\bbZ$. However, the following question remains relatively wide open.

\begin{qtn}\label{qtn:PU?}
Let $\Gam$ be a lattice in the unique cover $G_d$ of $\PU(n,1)$ of degree $1 \le d \le \infty$. Is $\Gam$ residually finite?
\end{qtn}

The cases $d = 1$ and $d = n+1$ are certainly yes: $\PU(n,1)$ and its $(n+1)$-fold cover $\SU(n,1)$ are linear. Proposition \ref{prop:FiniteCoverRF} implies that the same holds for the cover of $\PU(n,1)$ of any degree dividing $n + 1$, since it is intermediate to the covering by $\SU(n,1)$. Much less trivially, Question \ref{qtn:PU?} has a positive answer for the best-understood infinite class of commensurability classes of lattices in $\PU(n,1)$.

\begin{thm}[\cite{StoverToledo2, Hill}]\label{thm:PURF}
Suppose that $n \ge 2$ and $\Gam < \PU(n,1)$ is a cocompact lattice that is arithmetic of simplest kind. Then the preimage of $\Gam$ in any connected cover of $\PU(n,1)$ is residually finite.
\end{thm}

Lattices of simplest kind are those for which the associated algebraic group is the special unitary group of a hermitian form over a number field. The methods from joint work with Domingo Toledo \cite{StoverToledo2} and those of Richard Hill \cite{Hill} are similar in the broad approach but are ultimately different. In a very precise sense, the proof in \cite{StoverToledo2} is a generalization of the argument given above for Theorem \ref{thm:SL2RF}. Indeed, Theorem \ref{thm:PURF} is an immediate consequence of Lemma \ref{lem:2StepQuo} and the following result.

\begin{thm}[\cite{StoverToledo2}]\label{thm:PUnil}
Suppose that $n \ge 2$ and $\Gam < \PU(n,1)$ is a cocompact arithmetic lattice of simplest kind. If $\wt{\Gam}$ is the preimage of $\Gam$ in the universal cover $\wt{G}$ of $\PU(n,1)$, then there is a two-step nilpotent group $\calN$ and a finite index subgroup $\wt{\Gam}_0$ of $\wt{\Gam}$ that admits a surjection $\wt{\Gam}_0 \to \calN$ that is injective on the center of $\wt{\Gam}_0$.
\end{thm}

The assumption that $\Gam$ is cocompact of the simplest kind is used to produce the nilpotent group $\calN$ using the structure of its cohomology ring. The first piece of the puzzle is the following general result.

\begin{thm}[Thm.\ 1.2 \cite{StoverToledo1}]\label{thm:KahlerCondition}
Let $\Gam < \PU(n,1)$ be a torsion-free cocompact lattice and $X = \Gam \bs \bbB^n$. If $\wt{\Gam}$ is the preimage of $\Gam$ in the universal cover of $\PU(n,1)$, then there is a two-step nilpotent quotient of $\wt{\Gam}$ that is injective on the center of $\wt{\Gam}$ if and only if the K\"ahler class on $X$ is in the image of the cup product map
\[
c_\bbQ : \bigwedge\nolimits^2 H^1(X, \bbQ) \lra H^2(X, \bbQ).
\]
\end{thm}

The connection between the K\"ahler class and the preimage in the universal cover is that the first Chern class $c_1(X)$ is proportional to the K\"ahler class and is also the class in $H^2(X, \bbZ)$ associated with the preimage in the universal cover. This is because $\wh{\Gam} \bs \SU(n,1) / \SU(n)$ is homotopy equivalent to the canonical bundle minus its zero section, where $\wh{\Gam}$ is the preimage in $\SU(n,1)$. See \cite[\S 3.4]{StoverToledo1} for details. This also completes the connection with the above proof of Theorem \ref{thm:SL2RF}, since the canonical bundle of a Riemann surface is the cotangent bundle. The connection between a cohomology class being in the image of the cup product and nilpotent quotients is through the following result from algebraic topology.

\begin{thm}[Thm.\ 5.1 \cite{StoverToledo1}]\label{thm:NilCup}
Let $X$ be a closed aspherical manifold with fundamental group $\Gam$ and
\[
\U(1) \lra Y \lra X
\]
be a principal $\U(1)$ bundle with Euler class $\om \in H^2(X, \bbZ)$ that has infinite additive order, equivalently, its image in $H^2(X, \bbQ)$ is not zero. If $\wt{\Gam} = \pi_1(Y)$ and $z$ denotes a generator for $\pi_1(\U(1)) < \wt{\Gam}$, then the image of $z$ in the maximal two-step nilpotent quotient of $\wt{\Gam}$ has infinite order if and only if $\om$ is in the image of the map
\[
c_\bbQ : \bigwedge\nolimits^2 H^1(X, \bbQ) \lra H^2(X, \bbQ)
\]
given by evaluation of the cup product.
\end{thm}

\begin{pf}[Sketch of the direction used in this paper]
Taking the $K(\Gam, 1)$ map from the natural homomorphism $\Gam \to H_1(\Gam, \bbQ)$ gives a $\pi_1$-surjective map from $X$ to a $d$-dimensional torus $T$. The hypothesis on the cup product precisely implies that $\om$ is the pull-back of a class in $H^2(T, \bbQ)$. The associated circle bundle $W$ over $T$ has two-step nilpotent fundamental group, and the induced map $Y \to W$ produces the required homomorphism.
\end{pf}

Note that Theorems \ref{thm:KahlerCondition} and \ref{thm:NilCup} are very general, and have nothing to do with arithmeticity. However, the assumption that the lattice is arithmetic is simplest kind is used in a fundamental way in showing that these results apply to prove Theorem \ref{thm:PUnil} and thus Theorem \ref{thm:PURF}. The key point is that one can understand $H^1$ on congruence subgroups via the \emph{theta correspondence}; see \cite[\S 3]{StoverToledo2} for definitions and discussion. Closely following previous work of Bergeron, Li, Millson, and Moeglin \cite[Cor.\ 8.4]{BLMM} that applies Kudla--Millson theory (e.g., \cite{KudlaMillson}) in a clever way, joint work with Toledo proved the following first step that relates the class in $H^2(\Gam, \bbZ)$ associated with the canonical bundle with the geometry of the space $\Gam \bs \bbB^n$.

\begin{thm}[Thm.\ 3.3 \cite{StoverToledo2}]\label{thm:ChernInSpan}
Let $X = \Gam \bs \bbB^n$ be a smooth compact ball quotient, $n \ge 2$, with $\Gam$ a congruence arithmetic lattice of simple type, and let $c_1$ be the Chern form on $X$, considered as an element of $H^2(X, \bbC)$. Then $c_1$ is contained in the span of the collection of all Poincar\'e duals to codimension one totally geodesic subvarieties of $X$.
\end{thm}

In short, Theorem \ref{thm:ChernInSpan} says that the class in $H^2(\Gam \bs \bbB^n, \bbC)$ associated with its preimage in the universal cover of $\PU(n,1)$ will (virtually) be in the image of the cup product from $H^1(\Gam \bs \bbB^n, \bbC)$ if the same is true for any given finite collection of Poincar\'e duals to immersed totally geodesic subvarieties of (complex) codimension one. This is precisely the sort of result one obtains from the theta correspondence. Now following work of Bergeron, Millson, and Moeglin, joint work with Toledo proved the following result that, following the previous outline given in this section, proves Theorem \ref{thm:PURF}.

\begin{thm}[Thm.\ 3.4 \cite{StoverToledo2}]\label{thm:Cup}
Let $X = \Gam \bs \bbB^n$ be a compact congruence arithmetic ball quotient of simple type and $\mathrm{SC}^1(X) \subseteq H^2(X, \bbC)$ be the span of the Poincar\'e duals to the totally geodesic codimension one subvarieties of $X$. Then:
\begin{enum}

\item For any class $\sig \in \mathrm{SC}^1(X)$, there is a congruence cover $p : X^\prime \to X$ so that $p^*(\sig) \in \mathrm{SC}^1(X^\prime)$ is in the image of the cup product map
\[
c : \bigwedge\nolimits^2 H^1(X^\prime, \bbC) \lra H^2(X^\prime, \bbC).
\]

\item If $n \ge 3$, then for every $\phi \in H^{1,1}(X, \bbC)$ there is a congruence cover $p : X^\prime \to X$ so that $p^*(\phi)$ is contained in the image of $c$.

\item If $n \ge 4$, then for all $\phi \in H^2(X, \bbC)$ there is a congruence cover $p : X^\prime \to X$ so that $p^*(\phi)$ is in the image of $c$.

\end{enum}

\end{thm}

\begin{rem}\label{rem:StoverToledoMore}
Theorem \ref{thm:Cup} and the methods described in this section lead to residual finiteness of many other central extensions of arithmetic lattices in $\PU(n,1)$ of simplest kind. See \cite[\S 1]{StoverToledo2} for further results and applications, including the first construction of higher-dimensional smooth projective varieties that admit a K\"ahler metric of strongly negative curvature and are not homotopy equivalent to any locally symmetric space.
\end{rem}

\subsubsection{Remaining cases for $\PU(n,1)$}

The remaining cases for $\PU(n,1)$ are best subdivided into
\begin{itemize}

\item[$\star$] nonuniform arithmetic lattices,

\item[$\star$] arithmetic lattices constructed using hermitian forms on noncommutaive division algebras, and

\item[$\star$] nonarithmetic lattices.

\end{itemize}
In each case there are technical difficulties of varying significance in adapting the techniques described in this section. The only known result is for a single commensurability class.

\begin{thm}[\cite{StoverToledo1}, \cite{Hill}]\label{thm:Eisenstein}
Let $\Gam < \PU(2,1)$ be any lattice commensurable with the Picard modular group over the field $\bbQ(\sqrt{-3})$ and $\wt{\Gam}$ be its preimage in the universal cover of $\PU(2,1)$. Then $\wt{\Gam}$ is residually finite and linear.
\end{thm}

The proof given in joint work with Toledo is a computation much like the above proof of Theorem \ref{thm:SL2RF}, and does imply residual finiteness of the preimage in any connected cover of $\PU(2,1)$. See \cite[\S 8]{StoverToledo1}. Hill's proof is also computational, though quite different in its details. It seems reasonable to expect that the preimage of any nonuniform arithmetic subgroup of $\PU(n,1)$ in any of its connected covers will be residually finite.

It is much more difficult to predict what will be the case for the remainder of the arithmetic lattices, which are constructed using hermitian forms on vector spaces over noncommutaive division algebras. Indeed, the above strategy falters at the first step, since $H^1(\Gam, \bbZ)$ is trivial for any congruence subgroup by work of Rogawski \cite[Thm.\ 15.3.1]{Rogawski}, Clozel \cite[\S 3]{Clozel}, and Rajan--Venkataramana \cite[Thm.\ 4]{RajanVenky}. If preimages of these lattices in covers of $\SU(n,1)$ are residually finite, a new idea will be required, either by proving that there are noncongruence subgroups with the requisite cohomology to use the cup product strategy of \cite{StoverToledo1, StoverToledo2} (which would itself be profoundly interesting) or by using completely different methods. It is perhaps worth singling out one special case of interest that comes from this arithmetic construction \cite{PrasadYeung, PrasadYeungAdd, CartwrightSteger}.

\begin{qtn}\label{qtn:FakeP2}
Let $X$ be a fake projective plane with fundamental group $\Gam < \PU(2,1)$. Is the preimage of $\Gam$ in the every connected covering of $\PU(2,1)$ residually finite?
\end{qtn}

The case of the universal cover of $\PU(2,1)$ is particularly interesting, since the preimage $\wt{\Gam}$ of the fundamental group of the fake projective plane is then the fundamental group of the canonical line bundle minus its zero section. This space is homotopy equivalent to a circle bundle over the fake projective plane, and a relatively easy calculation with the associated Gysin sequence implies that $\wt{\Gam}$ is the fundamental group of a rational homology $5$-sphere. The existence of a homology sphere with fundamental group that is not residually finite would be quite interesting in and of itself. Since there are known presentations for these lattices by work of Cartwright and Steger \cite{CartwrightSteger}, it is possible that Question \ref{qtn:FakeP2} could be answered by computational methods.

\medskip

Moving on to nonarithmetic lattices, again apparently nothing is known. There are finitely many known commensurability classes; see \cite{DPP, Deraux3d} for a complete list of the known examples.  Presentations are known for the existing examples, so one could perhaps prove residual finiteness theorems using computational methods analogous to the arguments in \cite[\S 8]{StoverToledo1}. To ask a specific question that seems tractable:

\begin{qtn}\label{qtn:DM}
Let $\Gam < \PU(n,1)$ be one of the nonarithmetic lattices constructed by Deligne and Mostow \cite{DeligneMostow, MostowINT} and $\wt{\Gam}$ be its preimage in the universal cover of $\PU(n,1)$. Is there a finite index subgroup $\wt{\Lam}$ of $\wt{\Gam}$ and a nilpotent quotient $\calN$ of $\wt{\Lam}$ that is injective on the center?
\end{qtn}

These lattices tend to have finite index subgroups with infinite abelianization, so the cup product methods used in \cite{StoverToledo1, StoverToledo2} may be successful in answering Question \ref{qtn:DM}. What is much less clear, and would be independently interesting to know, is whether the K\"ahler class is in the span of the totally geodesic divisors on a nonarithmetic ball quotient. This ends up being quite subtle since, unlike the arithmetic setting, there are only finitely many totally geodesic divisors on nonarithmetic ball quotients \cite[Cor.\ 1.2]{BFMS2}.

\subsubsection{The wide open case of $\mathrm{F}_4^{(-20)}$}

One final case remains in rank one, and it is quite mysterious. Since the universal cover of $\mathrm{F}_4^{(-20)}$ is a nonlinear (double) cover, it is unclear whether or not lattices there are residually finite.

\begin{qtn}\label{qtn:F4}
Let $G$ be the universal cover of the adjoint real Lie group $\mathrm{F}_4^{(-20)}$ of type $\mathrm{F}_4$ and $\Gam < G$ be a lattice. Is $\Gam$ residually finite?
\end{qtn}

Lattices as in Question \ref{qtn:F4} admit a central exact sequence
\[
1 \lra L \lra \Gam \lra \Gam_0 \lra 1
\]
with $L \le \bbZ / 2$ and $\Gam_0$ a lattice in the linear group $\mathrm{F}_4^{(-20)}$. Therefore, by Corollary \ref{cor:JustCenter}, it suffices to prove that there is a finite quotient of $\Gam$ that is injective on its center. It is not clear how one would find such a quotient, and seems difficult at the moment to predict whether the answer to Question \ref{qtn:F4} should be yes or no.

\subsubsection{Gromov hyperbolic groups}

One of the major open problems in geometric group theory is whether or not every Gromov hyperbolic group is residually finite; see for example \cite[Prob.\ 3, p.\ 208]{GGT}. It is widely expected that there are hyperbolic groups that are not residually finite, but an example has yet to be discovered. The direct connection with this survey is that cocompact rank one lattices are Gromov hyperbolic groups. Indeed, they are virtually the fundamental group of compact manifolds admitting a Riemannian metric of negative curvature, which are the quintessential hyperbolic groups. However, something slightly stronger is true thanks to relatively hyperbolic Dehn filling. The following is well-known to experts.

\begin{prop}\label{prop:Gromov}
Let $\Gam$ be a lattice in a connected rank one real Lie group $G$ with finite center. If all Gromov hyperbolic groups are residually finite, then $\Gam$ is residually finite. If $\Gam$ is not residually finite, then there is a Gromov hyperbolic group that is not residually finite.
\end{prop}

\begin{pf}
If $\Gam$ is cocompact, then both conclusions of the proposition are immediate, since $\Gam$ is itself hyperbolic. Therefore, assume that $\Gam$ is nonuniform. Then there is a finite collection of peripheral subgroups $\calP_i$ in $\Gam$ so that $\Gam$ is hyperbolic relative to the collection $\{\calP_i\}$ \cite[Thm.\ 5.1]{Farb}. Note that $\Gam$ fits into a central extension
\[
1 \lra Z \lra \Gam \lra \overline{\Gam} \lra 1,
\]
where $\overline{\Gam}$ is a lattice in the adjoint form of $G$ and $Z$ is finite. Since $\overline{\Gam}$ is linear and thus residually finite, residual finiteness of $\Gam$ is equivalent to the existence of a finite quotient of $\Gam$ that is injective on $Z$ by Corollary \ref{cor:JustCenter}.

The groups $\calP_i$ are virtually nilpotent, so there are finite index normal subgroups $\calQ_i \le \calP_i$ that intersect $Z$ trivially. Choosing the $\calQ_i$ sufficiently deep inside $\calP_i$, a theorem of Osin \cite[Thm.\ 1.1]{Osin} implies that the quotient
\[
\Lam = \Gam / \ll\! \{\calQ_i\}\! \gg
\]
is hyperbolic relative to the finite groups $\{\calP_i / \calQ_i\}$ and $\Gam \to \Lam$ is injective on $Z$. Since finite groups are hyperbolic, $\Lam$ is a hyperbolic group \cite[Cor.\ 1.2]{Osin}.

If every hyperbolic group, in particular $\Lam$, is residually finite, then there is a finite quotient of $\Lam$ that is injective on $Z$ by Lemma \ref{lem:FiniteInject}. As mentioned earlier in the proof, this implies that $\Gam$ is residually finite. Conversely, if $\Gam$ is not residually finite, then there is no finite quotient of $\Gam$ that is injective on $Z$. Since $\Gam \to \Lam$ is injective on $Z$, then $\Lam$ also has no such quotient and therefore it is a hyperbolic group that is not residually finite. This proves the proposition.
\end{pf}

\begin{rem}
One can also extend Proposition \ref{prop:Gromov} to include the last remaining case where the center is infinite, namely the universal cover $\wt{G}$ of $\PU(n,1)$ for $n \ge 2$. Note that even cocompact lattices in $\wt{G}$ fail to be Gromov hyperbolic, since they have infinite center. The key observation is simply that if $\wt{\Gam}$ is a lattice in $\wt{G}$ and it is not residually finite, then its projection to some finite cover of $\PU(n,1)$ is also not residually finite.
\end{rem}

\section{The general case}\label{sec:Gen}

With all the above information on lattices in solvable and semisimple groups in hand, this section is actually quite short. Suppose that $G$ is a connected real Lie group and $\Lam < G$ is a lattice. If $G$ is simply connected, then it has a \emph{Levi decomposition}
\[
G = S \rtimes H
\]
where $S$ and $H$ are simply connected, $S$ is solvable, and $H$ is semisimple \cite[\S 1.4]{RussianLie}. When $G$ is not simply connected, one only obtains an exact sequence with connected solvable kernel and semisimple quotient.

As seen throughout this article, it generally suffices to assume that $G$ is simply connected and study quotients only when necessary. For simplicity it is a therefore standing assumption throughout this section that $G$ is simply connected. The first item of business is then to connect the Levi decomposition of $G$ to $\Lam$, which is subtle in general. The following theorem gives a sufficient condition for $\Lam$ to be the semidirect product of a lattice in $H$ with a lattice in $S$.

\begin{thm}\label{thm:FindSemidirect}
Suppose that $G$ is a simply connected Lie group with Levi decomposition $S \rtimes H$, and assume that $H$ has no compact factor. Then any lattice $\Lam < G$ can be decomposed as a semidirect product $\Del \rtimes \Gam$, where $\Del$ is a lattice in $S$ and $\Gam$ is a lattice in $H$.
\end{thm}

See \cite{Geng} for history and subtleties, including discussion of some flaws in \cite[Cor.\ 8.28]{RaghunathanBook} related only to the case where $H$ has compact factors. Even for general $G$, where the semisimple part is technically a quotient and may have compact factors, $\Lam$ still intersects $S$ in a lattice if and only if its projection to $H$ is a lattice; see \cite[Thm.\ 2.6]{Geng} for a reference. This leads to the following very general positive result.

\begin{thm}\label{thm:GenRF}
Let $G$ be a simply connected Lie group with Levi decomposition $S \rtimes H$ such that the semisimple part $H$ has no compact factor. If every lattice in $H$ is residually finite, then every lattice in $G$ is residually finite.
\end{thm}

\begin{pf}
If $\Lam < G$ is a lattice, then Theorem \ref{thm:FindSemidirect} implies that $\Lam = \Del \rtimes \Gam$ is a semidirect product of a lattice in $S$ and a lattice in $H$. Since $\Del$ is residually finite by Theorem \ref{thm:LinearSolvable} and finitely generated, Lemma \ref{lem:SemiDirect} implies that $\Lam$ is residually finite if and only if $\Gam$ is residually finite. This proves the theorem.
\end{pf}

\begin{rem}
Note that residual finiteness of all lattices in $G$ does not imply that every lattice in $H$ is residually finite. This conclusion would require knowing that for every lattice $\Gam_0 < H$, there is a finite index subgroup $\Gam \le \Gam_0$ and a lattice $\Del < S$ so that there is a lattice in $G$ isomorphic to $\Del \rtimes \Gam$. This statement often fails badly. For a well-known example, suppose $\Gam$ is a cocompact lattice in $\SL_n(\bbR)$ for $n \ge 3$. Then there is no lattice in ${\bbR^n \rtimes \SL_n(\bbR)}$ associated with $\Gam$ in this fashion. Indeed, such a lattice would necessarily be isomorphic to $\bbZ^n \rtimes \Gam$, and thus $\Gam$ is conjugate to a finite-index subgroup of $\SL_n(\bbZ)$, which is impossible.
\end{rem}

For an easy example of a lattice that is not residually finite, consider the semidirect product
\[
\bbZ^{2 g} \rtimes \wt{\Gam} < \bbR^{2 g} \rtimes \wt{G},
\]
where $\wt{\Gam}$ is the preimage of $\Sp_{2g}(\bbZ)$ in the universal cover $\wt{G}$ of $\Sp_{2g}(\bbR)$ and $\wt{G}$ acts on $\bbR^{2 g}$ through the standard action of $\Sp_{2 g}(\bbR)$. Similar examples lead to a wealth of lattices in general Lie groups that are not residually finite, but for fairly trivial reasons given the semisimple case. Indeed, perhaps the main takeaway of this section is that the semisimple part of the Levi decomposition is the primary deciding factor in whether a lattice in $G$ is residually finite\footnote{This is equally the case when $H$ has compact factors, where the projection of $\Lam$ to $H$ may not be discrete and thus has a slim chance to not be residually finite even if every lattice in $H$ is residually finite. It would be quite surprising if such an example exists.}. Hopefully, with the previous results of this paper in hand the reader is now well-prepared to enumerate the numerous other possibilities, theorems, examples, and pathologies regarding residual finiteness of lattices in general Lie groups.

\section*{Appendix: Fundamental groups of simple Lie groups}

This appendix tabulates fundamental groups of almost simple real Lie groups. Table \ref{tb:Compactpi1} provides the fundamental groups of various compact Lie groups $K$, the center $Z(\wt{K})$ of its universal cover $\wt{K}$ (i.e., the largest fundamental group of a group locally isomorphic to $\wt{K}$), and the index of $\pi_1(K)$ in $Z(\wt{K})$. Tables \ref{tb:Classicalpi1} and \ref{tb:Exceptionalpi1} use that data and Theorem \ref{thm:Cartan} to list the fundamental groups of classical and exceptional absolutely almost simple, absolutely simply connected linear algebraic groups $G$ and their adjoint form $\conj{G}$, along with applying Lemma \ref{lem:UCLinear} to record whether or not the universal cover $\wt{G}$ is linear.

\renewcommand{\arraystretch}{1.3}

\begin{table}[h]
\centering
\begin{tabular}{|c|c|c|}
\hline
$K$ & $Z(\wt{K})$ & $[Z(\wt{K}) : \pi_1(K)]$ \\
\hline
$\SU(n)$, $n \ge 2$ & $\bbZ / n$ & $1$ \\
\hline
$\U(n)$, $n \ge 1$ & $\bbZ$ & $n$ \\
\hline
$\SO(2)$ & $\bbZ$ & $2$ \\
\hline
$\SO(2 n + 1),\ n \ge 1$ & $\bbZ / 2$ & $1$ \\
\hline
$\SO(2 n),\ n \ge 2$ even & $(\bbZ / 2)^2$ & $2$ \\
\hline
$\SO(2 n),\ n \ge 1$ odd & $\bbZ / 4$ & $2$ \\
\hline
$\Sp(n)$ & $\bbZ / 2$ & $2$ \\
\hline
$\mathrm{E}_6$ & $\bbZ / 3$ & $3$ \\
\hline
$\mathrm{E}_7$ & $\bbZ / 2$ & $2$ \\
\hline
$\mathrm{E}_8$ & $\{1\}$ & $1$ \\
\hline
$\mathrm{F}_4$ & $\{1\}$ & $1$ \\
\hline
$\mathrm{G}_2$ & $\{1\}$ & $1$ \\
\hline
\end{tabular}
\caption{Compact Lie groups and their fundamental groups}\label{tb:Compactpi1}
\end{table}

\begin{table}[h]
\centering
\scriptsize{
\begin{tabular}{|c|c|c|c|}
\hline
$G$ & $K$ & $\pi_1(\conj{G})$ & $\wt{G}$ linear \\
\hline
$\SL_2(\bbR)$ & $\SO(2)$ & $\bbZ$ & no \\
\hline
$\SL_{2 n + 1}(\bbR)$, $n \ge 1$ & $\SO(2n+1)$ & $\bbZ / 2$ & no \\
\hline
$\SL_{2 n}(\bbR)$, $n \ge 2$ even & $\SO(2n)$ & $(\bbZ / 2)^2$ & no \\
\hline
$\SL_{2 n}(\bbR)$, $n \ge 2$ odd & $\SO(2n)$ & $\bbZ / 4$ & no \\
\hline
$\SL_n(\bbH)$, $n \ge 2$ & $\Sp(n)$ & $\bbZ / 2$ & yes \\
\hline
$\SU(p,q)$, $p \ge q \ge 1$ & $\mathrm{S}(\U(p) \times \U(q))$ & $\bbZ$ & no \\
\hline
$\SO(2 n + 1, 2)^\circ$, $n \ge 1$ & $\SO(2 n + 1) \times \SO(2)$ & $\bbZ/2 \times \bbZ$ & no \\
\hline
$\SO(1, 2 q)^\circ$, $q \ge 2$ even & $\SO(2 q)$ & $(\bbZ / 2)^2$ & yes \\
\hline
$\SO(2 p + 1, 2 q)^\circ$, $p \ge 0$, $q \ge 2$ even & $\SO(2 p + 1) \times \SO(2 q)$ & $(\bbZ / 2)^3$ & no \\
\hline
$\SO(1, 2 q)^\circ$, $q \ge 2$ odd & $\SO(2 q)$ & $\bbZ / 4$ & yes \\
\hline
$\SO(2 p + 1, 2 q)^\circ$, $p \ge 0$, $q \ge 2$ odd & $\SO(2 p + 1) \times \SO(2 q)$ & $\bbZ / 2 \times \bbZ / 4$ & no \\
\hline
$\Sp_{2n}(\bbR)$, $n \ge 2$ & $\U(n)$ & $\bbZ$ & no \\
\hline
$\Sp(p, q)$, $p \ge q \ge 1$ & $\Sp(p) \times \Sp(q)$ & $(\bbZ / 2)^2$ & yes \\
\hline
$\SO(2 p + 1, 1)^\circ$, $p > 0$ & $\SO(2 p + 1)$ & $\bbZ / 2$ & yes \\
\hline
$\SO(2 p + 1, 2 q + 1)^\circ$, $p \ge q \ge 0$, $p > 0$ & $\SO(2 p + 1) \times \SO(2 q + 1)$ & $(\bbZ / 2)^2$ & no \\
\hline
$\SO(2 n, 2)^\circ$, $n \ge 2$ even & $\SO(2 n) \times \SO(2)$ & $(\bbZ / 2)^2 \times \bbZ$ & no \\
\hline
$\SO(2 n, 2)^\circ$, $n \ge 3$ odd & $\SO(2 n) \times \SO(2)$ & $\bbZ / 4 \times \bbZ$ & no \\
\hline
$\SO(2 p, 2 q)^\circ$, $p \ge q \ge 2$ even & $\SO(2 p) \times \SO(2 q)$ & $(\bbZ / 2)^4$ & no \\
\hline
$\SO(2 p, 2 q)^\circ$, $p, q \ge 2$, $p$ odd, $q$ odd & $\SO(2 p) \times \SO(2 q)$ & $(\bbZ / 4)^2$ & no \\
\hline
$\SO(2 p, 2 q)^\circ$, $p, q \ge 2$, $p$ even, $q$ odd & $\SO(2 p) \times \SO(2 q)$ & $(\bbZ / 2)^2 \times \bbZ / 4$ & no \\
\hline
$\SO^*(2 n)$, $n \ge 2$ & $\U(n)$ & $\bbZ$ & no \\
\hline
\end{tabular}
}
\caption{Classical simple Lie groups}\label{tb:Classicalpi1}
\end{table}

\begin{table}[h]
\centering
\begin{tabular}{|c|c|c|c|}
\hline
$G$ & $K$ & $\pi_1(\conj{G})$ & $\wt{G}$ linear \\
\hline
$\mathrm{E}_6^{(6)}$ & $\Sp(4)$ & $\bbZ / 2$ & no \\
\hline
$\mathrm{E}_6^{(2)}$ & $\SU(6) \times \SU(2)$ & $\bbZ / 6$ & no \\
\hline
$\mathrm{E}_6^{(-14)}$ & $\mathrm{Spin}(10) \times \SO(2)$ & $\bbZ$ & no \\
\hline
$\mathrm{E}_6^{(-26)}$ & $\mathrm{F}_4$ & $\{1\}$ & yes \\
\hline
$\mathrm{E}_7^{(7)}$ & $\SU(8)$ & $\bbZ / 4$ & no \\
\hline
$\mathrm{E}_7^{(-5)}$ & $\SU(2) \times \SO(12)$ & $(\bbZ / 2)^2$ & no \\
\hline
$\mathrm{E}_7^{(-25)}$ & $\SO(2) \times \mathrm{E}_6$ & $\bbZ$ & no \\
\hline
$\mathrm{E}_8^{(8)}$ & $\mathrm{Spin}(16)$ & $\bbZ / 2$ & no \\
\hline
$\mathrm{E}_8^{(-24)}$ & $\SU(2) \times \mathrm{E}_7$ & $\bbZ / 2$ & no \\
\hline
$\mathrm{F}_4^{(4)}$ & $\SU(2) \times \Sp(3)$ & $\bbZ / 2$ & no \\
\hline
$\mathrm{F}_4^{(-20)}$ & $\mathrm{Spin}(9)$ & $\bbZ / 2$ & no \\
\hline
$\mathrm{G}_2^{(2)}$ & $\SU(2) \times \SU(2)$ & $\bbZ / 2$ & no \\
\hline
\end{tabular}
\caption{Exceptional simple Lie groups}\label{tb:Exceptionalpi1}
\end{table}

\renewcommand{\arraystretch}{1}

\clearpage

\bibliography{RFSurvey}

\end{document}